\documentclass[a4paper,12pt,leqno]{amsart} 
\usepackage{amsmath,amsthm,amssymb}  
\usepackage{mathdots}

\usepackage{graphicx}

\numberwithin{equation}{section}

\theoremstyle{plain}
\newtheorem{theorem}{Theorem}[section]
\newtheorem{proposition}[theorem]{Proposition}         
\newtheorem{corollary}[theorem]{Corollary}

\theoremstyle{definition}  
 
\newcommand{\C}{\mathbb C}   
\newcommand{\R}{\mathbb R}
\newcommand{\Z}{\mathbb Z}

\renewcommand{\P}{\mathbb P}
\newcommand{\X}{\mathbb X}  

\newcommand{\al}{\alpha}
 
\newcommand{\ga}{\gamma}
\newcommand{\de}{\delta}
\newcommand{\la}{\lambda}
\newcommand{\si}{\sigma} 
 
\newcommand{\La}{\Lambda}

\DeclareMathOperator{\diag}{diag}

\newcommand{\calM}{\mathcal{M}}

\newcommand{\no}{\noindent}

\newcommand{\pr}{\prime} 
\newcommand{\prr}{{\prime\prime}} 
\newcommand{\st}{\ \vert\ }   
\renewcommand{\ll}{\lq\lq}
\newcommand{\rr}{\rq\rq\ }
\newcommand{\rrr}{\rq\rq}  
\renewcommand{\b}{\partial}
\newcommand{\zbar}{  {\bar z}  }
\newcommand{\zzb}{ {z\bar z}  }

\newcommand{\bp}{\begin{pmatrix}} 
\newcommand{\ep}{\end{pmatrix}}

\newcommand{\GL}{\textrm{GL}}
\newcommand{\SL}{\textrm{SL}}
\renewcommand{\O}{\textrm{O}}

\newcommand{\GLNR}{\GL_n\R}

\newcommand{\sss}{\scriptsize}

\begin{document}     

\title[Integral Stokes data]{Some tt* structures\\ and their integral Stokes data}  

\author{Martin A. Guest and Chang-Shou Lin}          

\date{}   

\begin{abstract}  In \cite{GuItLiXX} 
a description was given of all smooth solutions 
of the two-function tt*-Toda equations in terms
of  asymptotic data,  holomorphic data, and 
monodromy data.  
In this supplementary article we focus on the holomorphic data
and its interpretation in quantum cohomology,
and enumerate those solutions with integral Stokes data.   This
leads to a characterization of quantum D-modules for certain
complete intersections of Fano type in weighted projective spaces.
\end{abstract}

\subjclass[2000]{Primary 81T40;
Secondary 53D45, 35J60}

\maketitle

\section{The tt*-Toda equations}\label{ttt}

The tt* (topological---anti-topological fusion) equations were introduced by
S.~Cecotti and C.~Vafa in their work on deformations of quantum field theories with N=2 supersymmetry (section 8 of \cite{CeVa91}, and also \cite{CeVa92},\cite{CeVa92a}).
This has led to the development of an area known as tt* geometry 
(\cite{CeVa91},\cite{Du93},\cite{He03}), a generalization of special geometry.

Solutions of the tt* equations can be interpreted as pluriharmonic maps with values in the noncompact real symmetric space $\GLNR/\O_n$, or as
pluriharmonic maps with values in a certain classifying space of  variations of polarized (finite or infinite-dimensional) Hodge structure.   Frobenius manifolds with real structure, e.g.\ quantum cohomology algebras, provide a very special class of solutions \ll of geometric origin\rr (see \cite{Du93}).   
These special solutions lie at the intersection of p.d.e.\  theory, integrable systems, and (differential, algebraic, and symplectic) geometry.  However, very few concrete examples have been worked out in detail, and their study is just beginning.   It is relatively straightforward to obtain local solutions, but these special solutions have (or are expected to have) global properties, and these properties are hard to establish.

In \cite{GuLiXX},  \cite{GuItLiXX}  a family of global solutions was constructed by relatively elementary p.d.e.\ methods. In this article we shall describe the special solutions in terms of their holomorphic data.   This allows us to obtain --- in a very restricted situation --- an a fortiori  characterization of quantum D-modules by purely algebraic/analytic means, which is one of the long term goals of the subject (cf.\ \cite{Gu08} and the extensive literature on o.d.e.\ of Calabi-Yau type).  

The equations studied in \cite{GuLiXX},  \cite{GuItLiXX} are
\begin{equation}\label{ost}
 2(w_i)_{\zzb}=-e^{2(w_{i+1}-w_{i})} + e^{2(w_{i}-w_{i-1})}
\end{equation}
where each $w_i$ is real-valued on (an open subset of) $\C=\R^2$, and $w_i=w_{i+n+1}$ for all $i\in\Z$; this is the two-dimensional periodic Toda lattice \ll with opposite sign\rrr.   In addition it is essential to assume  that
\begin{equation}\label{as}
\begin{cases}
\ \  w_0+w_{l-1}=0, \ w_1+w_{l-2}=0,\ \ \dots\\
\ \  w_l+w_{n}=0, \ w_{l+1}+w_{n-1}=0,\ \ \dots
 \end{cases}
\end{equation}
for some $l\in\{0,\dots,n+1\}$ (the cases $l=0$ and $l=n+1$ mean that
$w_i+w_{n-i}=0$ for all $i$).   The system (\ref{ost}), (\ref{as}) is then a special case of the tt* equations, and we call it the tt*-Toda system.  In the ten cases listed in Table \ref{tableofcases-short} below,  $w_0,\dots,w_n$ reduce to two unknown functions and  (\ref{ost}) reduces to
\begin{equation}\label{pde}
\begin{cases}
u_{z\zbar}&= \ e^{au} - e^{v-u} 
\\
v_{z\zbar}&= \ e^{v-u} - e^{-bv}
\end{cases}
\end{equation}
with $a,b\in\{1,2\}$, and it is this system that was solved in \cite{GuLiXX}
for $u,v:\C^\ast\to\R$.
\begin{table}[h]
\renewcommand{\arraystretch}{1.3}
\begin{tabular}{c||cc|cc|cc}\label{cases}
case & $l$ & $n\!\!+\!\!1\!\!-\!l$ & $u$ & $v$ & $a$ & $b$
\\
\hline
4a & $4$ & 0 & $2w_0$ & $2w_1$ & 2 & 2 
\\
4b & $2$ & $2$ & $2w_3$ & $2w_0$ & 2 & 2
\\
\hline
5a & $5$ & 0 & $2w_0$ & $2w_1$ & 2 & 1 
\\
5b & $3$ & $2$ & $2w_4$ & $2w_0$ & 2 & 1
\\
\hline
5c & $4$ & $1$ & $2w_0$ & $2w_1$ & 1 & 2 
\\
5d & $1$ & $4$ & $2w_1$ & $2w_2$ & 1 & 2  
\\
5e & $2$ & $3$ & $2w_4$ & $2w_0$ & 1 & 2 
\\
\hline
6a & $5$ & $1$ & $2w_0$ & $2w_1$ & 1 & 1 
\\
6b & $1$ & $5$ & $2w_1$ & $2w_2$ & 1 & 1 
\\
6c & $3$ & $3$ & $2w_5$ & $2w_0$ & 1 & 1
\end{tabular}
\bigskip
\caption{}\label{tableofcases-short}
\end{table}

\begin{figure}[h]
\begin{center}
\includegraphics[scale=0.6, trim= 0 300 0 250]{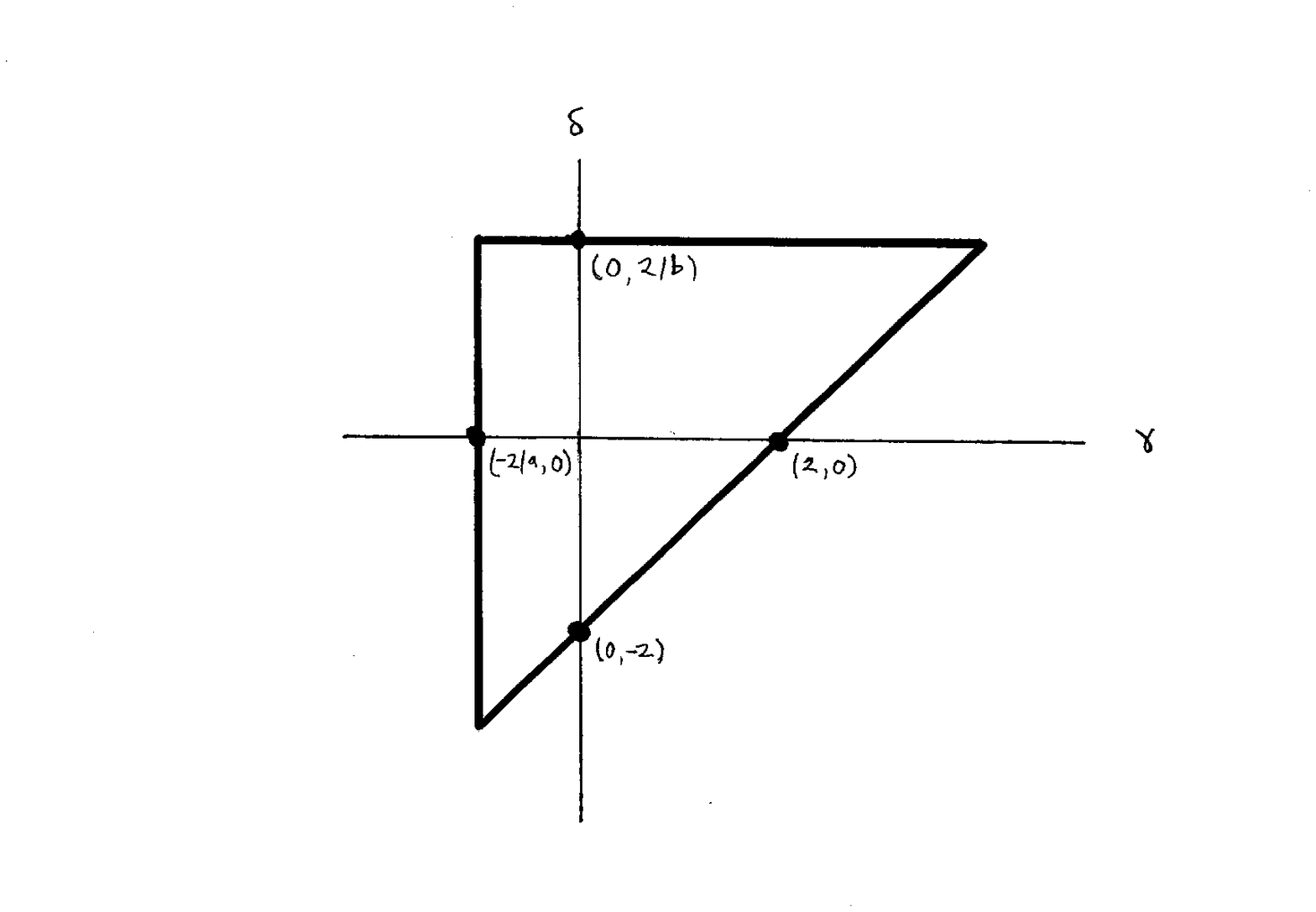}
\end{center}
\caption{}\label{theregion}
\end{figure}

The first main result of \cite{GuLiXX},  \cite{GuItLiXX} is a characterization of solutions of 
(\ref{pde}) in terms of asymptotic data.  Namely (Theorem A of \cite{GuItLiXX}),
for any $(\ga,\de)$ in the triangular region 
\[
\text{$\ga\ge -2/a$, $\de\le 2/b$, $\ga-\de\le2$}
\]
(Fig.\ \ref{theregion}),
the system
(\ref{pde})
has a unique solution $(u,v)$ such that
\begin{align*}
u(z) &\sim \ga\log \vert z\vert \\
v(z) &\sim \de \log \vert z\vert
\end{align*}
as $\vert z\vert\to 0$, and 
$u(z) \to 0$, $v(z) \to 0$ as $\vert z\vert\to \infty$.
The functions $u,v$ depend only on $\vert z\vert$. 

On the other hand, from the integrable systems point of view, these solutions correspond to two other kinds of data, which we shall describe next.

First, from the zero curvature formulation of (our version of) the Toda lattice, which may be written in the form 
\begin{align*}
F^{-1} F_z &=  \al^\pr\\
F^{-1} F_{\zbar} &=  \al^{\pr\pr},
\end{align*}
we have holomorphic data $\eta=L^{-1}L_z$ where $F=LB^{-1}$ is a Birkhoff factorization.
The system (\ref{ost}) is equivalent to $d\al+\al\wedge\al=0$, the condition that the connection $d+\al$ is flat. Here,
$\al=\al^\pr dz + \al^{\pr\pr} d\zbar$ is defined in terms of $w_0,\dots,w_n$.
We omit the details, which are given in \cite{GuLiXX}
and will be reviewed briefly in section \ref{qcoh}.

Next, for radial solutions, i.e.\ when the $w_i$ depend only on $x=\vert z\vert$, 
we can write $\al = \al^{\text rad} \, dx$, and
the flat connection $d+\al$ extends to a  flat connection $d+\al+\hat\al$, 
for some $\hat\al=\al^{\text sp} d\mu$, where $\mu$ is a \ll spectral parameter\rrr.
The radial version of system (\ref{ost}) is equivalent to the condition that the connection $d+\al+\hat\al$ is flat, and we obtain a rather different zero curvature formulation
\begin{align*}
F^{-1} F_x &=  \al^{\text rad} \\
F^{-1} F_{\mu} &=  \al^{\text sp}.
\end{align*}
Here $\al^{\text sp}$ is meromorphic in $\mu$ with poles of order $2$ at $\mu=0$ and 
$\mu=\infty$ (formulae can be found in \cite{GuItLiXX}).
The first equation can be regarded as describing an isomonodromic family of 
$x$-deformations of the second equation.    In particular the Stokes data is independent of $x$.  

It turns out that the Stokes data alone parametrizes the above solutions $u,v$ and in fact this Stokes data
reduces to two real numbers $s_1^\R,s_2^\R$.    The relation between the asymptotic data $\ga,\de$ and the Stokes data $s_1^\R,s_2^\R$ is as follows (Theorem B of \cite{GuItLiXX}):

\smallskip

\no (i) Cases 4a, 4b:

\no$\pm s_1^\R = 2\cos \tfrac\pi4 {\scriptstyle (\ga+1)} +  2\cos \tfrac\pi4 {\scriptstyle(\de+3)}$

\no$-s_2^\R = 2+4\cos \tfrac\pi4 {\scriptstyle(\ga+1)} \, \cos \tfrac\pi4 {\scriptstyle(\de+3)}$

\no (ii) Cases 5a, 5b:

\no\ \ $s_1^\R = 1+2\cos \tfrac\pi5 {\scriptstyle(\ga+6)}  + 2\cos \tfrac\pi5 {\scriptstyle(\de+8)}$

\no$-s_2^\R = 2
+2\cos \tfrac\pi5 {\scriptstyle(\ga+6)}  +  2\cos \tfrac\pi5 {\scriptstyle(\de+8)}
+4\cos \tfrac\pi5 {\scriptstyle(\ga+6)} \, \cos \tfrac\pi5 {\scriptstyle(\de+8)}$

\no (iii) Cases 5c, 5d, 5e:

\no
\ \ $s_1^\R = 1+2\cos \tfrac\pi5 {\scriptstyle(\ga+2)}  +  2\cos \tfrac\pi5 {\scriptstyle(\de+4)}$

\no$-s_2^\R = 2
+2\cos \tfrac\pi5 {\scriptstyle(\ga+2)}  +  2\cos \tfrac\pi5 {\scriptstyle(\de+4)}
+4\cos \tfrac\pi5 {\scriptstyle(\ga+2)} \, \cos \tfrac\pi5 {\scriptstyle(\de+4)}$

\no (iv) Cases 6a, 6b, 6c:

\no
$\pm s_1^\R = 2\cos \tfrac\pi6 {\scriptstyle(\ga+2)}  +  2\cos \tfrac\pi6 {\scriptstyle(\de+4)}$

\no$-s_2^\R = 1+4\cos \tfrac\pi6 {\scriptstyle(\ga+2)} \, \cos \tfrac\pi6 {\scriptstyle(\de+4)}$

\smallskip

The purpose of this article is to investigate the solutions for which $s_1^\R,s_2^\R$ are {\em integers.}   These are the \ll physical solutions\rr predicted by Cecotti and Vafa.  We shall describe them in terms of their holomorphic data, and (where appropriate) explain how this holomorphic data can be interpreted in terms of quantum cohomology or other geometrical phenomena.

In section \ref{qcoh} we review the definition of holomorphic data and compute it for the solutions of (\ref{pde}) described above.  In section \ref{integral} we identify those solutions with 
$s_1^\R,s_2^\R\in\Z$;  tables of all three types of data for all ten cases are presented in the appendix.   Quantum cohomology (or rather quantum D-module) interpretations of these solutions are given in section \ref{interp}.   This leads directly to the characterization result (Corollary \ref{char}).

\section{Holomorphic data for solutions}\label{qcoh}

The idea of holomorphic data for pluriharmonic maps has arisen independently in several contexts: 
in \cite{Kr80} as a method of solving integrable equations, in \cite{Se89} as a loop group version of the same thing,  in \cite{Si92} as a correspondence between harmonic bundles and holomorphic bundles,  and in \cite{DoPeWu98},\cite{DMPW97} as a systematic method for studying harmonic maps into symmetric spaces. And it appeared already in \cite{CeVa91} for the tt* equations themselves, though this was perhaps not appreciated at the time by differential geometers.

As a way of specifying a pluriharmonic map into a symmetric space, the holomorphic data generalizes the classical Weierstrass representation of a minimal surface.  It has the same advantages and the same disadvantages, and its usefulness depends on the circumstances.  In the case of the tt* equations, however, the holomorphic data plays a crucial  role, because of its
geometrical interpretation as a quantum cohomology ring of a manifold or Milnor ring of a singularity.

The holomorphic data
in our situation is (see \cite{GuLiXX}) a matrix of the form
\[
\eta=
\bp
 & & & p_0\\
 p_1 & & & \\
  & \ddots & & \\
   & & p_n &
   \ep
\]
where each $p_i=p_i(z)$ is a holomorphic function. 
From this holomorphic data we can construct {\em local} solutions
of (\ref{ost}) as follows.

For some $z_0\in U$ and some simply connected open neighbourhood $U^\prime$ of $z_0$ in $U$,  let $L$ be the solution of the holomorphic o.d.e.\ system $L^{-1}dL = \frac1\la\, \eta\, dz$, with initial condition $L(z_0)=I$.  We regard $L$ as a map $U^\pr\to \La\SL_{n+1}\C$, where $\La\SL_{n+1}\C$ is the free loop group of $\SL_{n+1}\C$, $\la$ being the loop parameter.
Let $L=FB$ be the
Iwasawa factorization of $L$ (see chapter 12 of \cite{Gu97}) with $F(z_0)=I$, $B(z_0)=I$. This factorization is possible on some neighbourhood $U^\prr$ of $z_0$. It follows that
$B$ is of the form $B=\sum_{i\ge 0}\la^i B_i$ where  $B_0=\diag(b_0,\dots,b_n)$. The factorization $L=FB$ is unique if we insist that $b_i>0$ for all $i$. We have $b_0\dots b_n=1$ and $b_i(z_0)=1$ for all $i$.

Let $\al = F^{-1}dF=F^{-1}F_z dz+F^{-1}F_{\zbar} d\zbar$. This
 must be of the form $\al^\pr dz + \al^\prr d\zbar$ where 
\[
\al^\pr=
\bp
a_0 & & & \\
 &  a_1 & & \\
  & & \ddots & \\
   & &  & a_n
\ep
+
\frac1\la
\bp
 & & & A_0\\
 A_1 & & & \\
  & \ddots & & \\
   & & A_n &
\ep
\]
for some smooth functions $a_i,A_j:U^\prr\to\C$.    
From the $\la^{-1}$ terms of $\al^\pr=F^{-1}F_z=(LB^{-1})^{-1}(LB^{-1})_z = \tfrac1\la B\eta B^{-1} + B(B^{-1})_z$, we obtain
$
A_i=p_i b_i/b_{i-1}
$
and similarly from the diagonal terms of $F^{-1}F_{\zbar}$ we obtain
$
a_i=(\log b_i)_z.
$
Since $\al = F^{-1}dF$, we have the zero curvature equation $d\al+ \al\wedge\al=0$, which gives an additional equation
\[
(a_i)_{\zbar} + (\bar a_i)_z
 = -\vert A_{i+1}\vert^2 + \vert A_{i}\vert^2.
\]

Let
$
w_i=\log b_i - \log \vert h_i\vert
$
where $h_0,\dots,h_n$ are any holomorphic functions. We obtain
\begin{equation}\label{provis}
2(w_i)_{\zzb}=
-\vert \nu_{i+1} \vert^2 e^{2(w_{i+1}-w_{i})} 
+ \vert \nu_{i} \vert^2 e^{2(w_{i}-w_{i-1})}
\end{equation}
where $\nu_i=p_i h_i/h_{i-1}$.  Choosing $h_0,\dots,h_n$ such that all $\nu_i$ are equal, say $\nu_i=\nu$ for all $i$, we have $\nu^{n+1}=p_0\dots p_n$ and $\nu=p_i h_i/h_{i-1}$.  

For consistency with (\ref{as}) we impose the condition that $h_ih_j=1$ whenever $w_i+w_j=0$ in (\ref{as}).  This determines $h_0,\dots,h_n$ explicitly in terms of $p_0,\dots,p_n$ (cf.\ section 4 of \cite{GuLiXX}).  

Finally, the change of variable $z\mapsto \int\nu\, dz$ then converts (\ref{provis})  into  (\ref{ost}).
We obtain the required local solution of (\ref{ost}), (\ref{as}).

Let us turn now to the holomorphic data for the solutions $w_i:\C^\ast\to\R$ 
parametrized by $(\ga,\de)$ in the triangular region 
of Fig.\ \ref{theregion}.   
The radial property implies that the holomorphic data must be of the form $p_i(z)=c_i z^{k_i}$.  

\begin{table}[h]
\renewcommand{\arraystretch}{1.3}
\begin{tabular}{c||c|c|l}
case & $N\ga$ & $N\de$ & $k_0,\dots,k_n$
 \\
\hline
$4a$  &  
${3k_0 - 2k_1 - k_2}$
&
${k_0 + 2k_1 - 3k_2}$
& 
$k_1=k_3$
\\
$4b$ &  
${-2k_0 - k_1 + 3k_3}$
&
${2k_0 - 3k_1 + k_3}$
& $k_0=k_2$
\\
\hline
$5a$  & ${4k_0 - 2k_1 -2k_2}$ 
&
${2k_0 + 4k_1 -6k_2}$
& 
$k_1=k_4, k_2=k_3$
\\
$5b$ &  
${-2k_0 - 2k_1 +4k_4}$
&
${4k_0 - 6k_1 +2k_4}$
& $k_0=k_3, k_1=k_2$
\\
\hline
$5c$ &  
${6k_0 -4k_1 -2k_2}$
&
${2k_0 +2k_1 -4k_2}$
& 
$k_1=k_3, k_0=k_4$
\\
$5d$&  
${6k_0 - 4k_2 - 2k_3}$
&
${2k_0 + 2k_2 - 4k_3}$
& 
$k_2=k_4, k_0=k_1$
\\
$5e$ &  
${-4k_0 -2k_1 +6k_3}$
&
${2k_0 -4k_1 +2k_3}$
& 
$k_0=k_2, k_3=k_4$
\\
\hline
$6a$ &  
${8k_0 -4k_1 -4k_2}$
&
${4k_0 +4k_1 -8k_2}$
& 
$k_1=k_4, k_0=k_5, k_2=k_3$
\\
$6b$ & 
${8k_0 -4k_2 -4k_3}$ 
&
${4k_0 +4k_2 -8k_3}$
& 
$k_2=k_5, k_0=k_1, k_3=k_4$
\\
$6c$ & 
${-4k_0 -4k_1 +8k_4}$ 
&
${4k_0 -8k_1 +4k_4}$
& 
$k_0=k_3, k_4=k_5, k_1=k_2$
\end{tabular}
\bigskip
\caption{}\label{holdata}
\end{table}

The relation between $k_0,\dots,k_n$ and $\ga,\de$
is given in Table \ref{holdata}.  This was obtained in section 4 of \cite{GuLiXX}, and it was explained there that one may normalize so that $c_0=\cdots=c_n=1$.
We write  $N=n+1+\sum_{i=0}^{n} k_i$ from now on.

Using this and Theorem A of \cite{GuItLiXX}, we find the following expressions for $s_1^\R,s_2^\R$
in terms of $k_0,\dots,k_n$.

\begin{proposition}\label{bprime} $ $

\no (i) Cases 4a ($k=k_0, l=k_2$), 4b ($k=k_3, l=k_1$): 

\no$\pm s_1^\R = 2\cos \tfrac\pi N {\scriptstyle (k+1)} -  2\cos \tfrac\pi N {\scriptstyle(l+1)}$

\no$-s_2^\R = 2-4\cos \tfrac\pi N {\scriptstyle(k+1)} \, \cos \tfrac\pi N {\scriptstyle(l+1)}$

\no (ii) Cases 5a ($k=k_0, l=k_2$), 5b ($k=k_4, l=k_1$):

\no\ \ $s_1^\R = 1-2\cos \tfrac\pi N {\scriptstyle(k+1)}  +  2\cos \tfrac{2\pi} N {\scriptstyle(l+1)}$

\no$-s_2^\R = 2
-2\cos \tfrac\pi N {\scriptstyle(k+1)}  +  2\cos \tfrac{2\pi} N {\scriptstyle(l+1)}
-4\cos \tfrac\pi N {\scriptstyle(k+1)} \, \cos \tfrac{2\pi} N {\scriptstyle(l+1)}$

\no (iii) Cases 5c ($k=k_0, l=k_2$), 5d ($k=k_0, l=k_3$), 5e ($k=k_3, l=k_1$):

\no
\ \ $s_1^\R = 1+2\cos \tfrac{2\pi}N {\scriptstyle(k+1)}  -  2\cos \tfrac{\pi} N {\scriptstyle(l+1)}$

\no$-s_2^\R = 2
+2\cos \tfrac{2\pi} N {\scriptstyle(k+1)}  -  2\cos \tfrac{\pi} N {\scriptstyle(l+1)}
-4\cos \tfrac{2\pi} N {\scriptstyle(k+1)} \, \cos \tfrac{\pi} N {\scriptstyle(l+1)}$

\no (iv) Cases 6a ($k=k_0, l=k_2$), 6b ($k=k_0, l=k_3$), 6c ($k=k_4, l=k_1$):

\no
$\pm s_1^\R = 2\cos \tfrac{2\pi} N {\scriptstyle(k+1)}  -  2\cos \tfrac{2\pi} N {\scriptstyle(l+1)}$

\no$-s_2^\R = 1-4\cos \tfrac{2\pi} N {\scriptstyle(k+1)} \, \cos \tfrac{2\pi} N {\scriptstyle(l+1)}$
\end{proposition}

These formulae make no reference to the \ll real structure\rrr, and in fact Proposition \ref{bprime} could have been obtained directly from the flat holomorphic connection $d+\tfrac1\la \eta dz$.  Namely, by homogeneity, $d+\tfrac1\la \eta dz$ extends to a flat connection $d+\tfrac1\la \eta dz + \hat\eta\,d\la$ (just as $d+\al$ extends to $d+\al+\hat\al$). The meromorphic connection $d+ \hat\eta\,d\la$ has poles of order $2,1$ at $\la=0,\infty$.  This is the connection usually considered in the theory of Frobenius manifolds.
The Stokes analysis of $d+ \hat\eta\,d\la$ at $\la=0$ is the same as that of $d+\hat\al$ at $\la=0$, because
\[
L=FB\sim F\ \ \text{as}\ \ \la\to 0.
\]
This leads to a relation between the Stokes data $s^\R_1,s^\R_2$ and the monodromy at the regular singular point.  The latter can be computed in terms of $k_0,\dots,k_n$, and the formulae above follow from this.  
This principle was used already in \cite{BoIt95} in a similar situation for $2\times 2$ matrices.   It shows that the Stokes matrices arising in the theory of Frobenius manifolds agree with those of the  tt* equations in the case where the Frobenius manifold admits a real structure.

However, the existence of a real structure --- equivalently, the existence of the Iwasawa factorization 
$L=FB$ --- is a nontrivial property.  It holds in our situation if and only if $k_0\ge-1,\dots,k_n\ge-1$.  This is a consequence of Theorem A of \cite{GuItLiXX}, as it may be deduced from Table \ref{holdata} that the conditions 
$\ga\ge -2/a, \de\le 2/b, \ga-\de\le2$ are equivalent to the conditions
$k_0\ge-1,\dots,k_n\ge-1$.

\section{Solutions with integral Stokes data}\label{integral}

\begin{figure}[h]
\begin{center}
\includegraphics[scale=0.6, trim= 0 300 0 240]{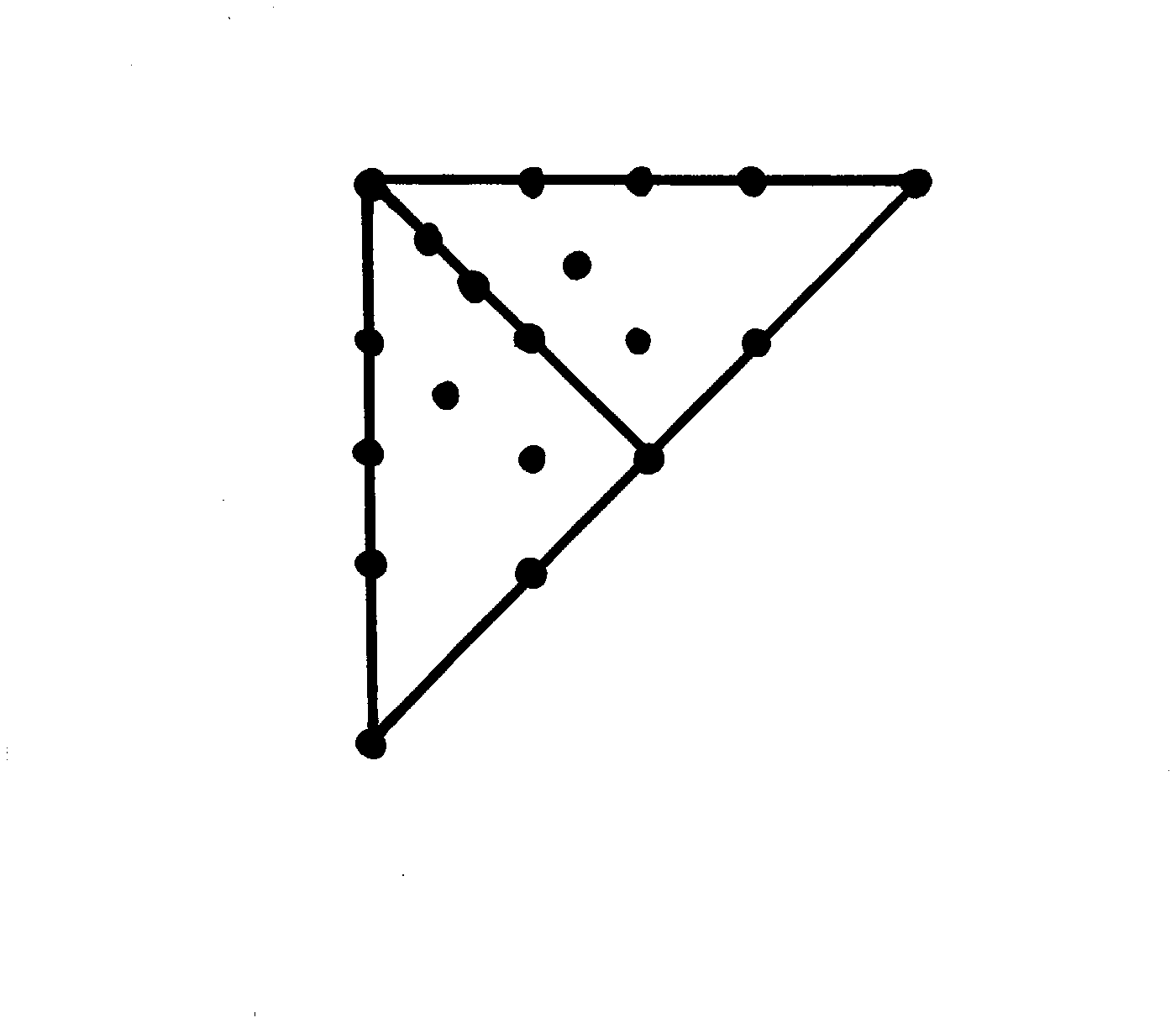}
\end{center}
\caption{The 19 points.}\label{thepoints}
\end{figure}

From the explicit formulae it is straightforward to to identify those solutions 
for which $s^\R_1,s^\R_2$ are integers:

\begin{proposition}
For each case in Table \ref{tableofcases-short}, 
there are 19 solutions $(u,v)$ for which the Stokes data $s_1^\R,s_2^\R$ is integral.
The corresponding values of the asymptotic data $(\ga,\de)$ are listed in Table \ref{integralsolutions}, and are shown schematically in Fig.\ \ref{thepoints}.  The Stokes data and holomorphic data for each of these solutions are given in Tables 
\ref{theoremCi}-\ref{theoremCiiii} of the appendix.
\end{proposition}
\begin{proof}   We use the formulae of Proposition \ref{bprime}.  
The region is given
by $k_i+1\ge0$ for all $i$ (see the end of section \ref{qcoh}). 
In all ten cases, the solutions with integral Stokes data are given by 
\[
2\cos a-2 \cos b\in\Z, \ \ 4 \cos a \cos b\in\Z
\]
where
$a=\tfrac\pi N {\scriptstyle(k+1)}$, $b=\tfrac\pi N {\scriptstyle(l+1)}$. An elementary
calculation shows that the set 
\[
\{(a,b)\in[0,\pi] \st 2\cos a-2\cos b\in\Z, 4\cos a\cos b\in\Z \}
\]
consists of the following 33 points: 25 points with $\cos a, \cos b\in \tfrac12\Z$, i.e.\ 
$a,b\in\{ 0, \tfrac\pi3, \tfrac\pi2, \tfrac{2\pi}3, \pi \}$, and 8 additional points
$(a,b)=(\tfrac{\pi}{6}, \tfrac{\pi}{6})$, $(\tfrac{\pi}{4}, \tfrac{\pi}{4})$, $(\tfrac{3\pi}{4}, \tfrac{3\pi}{4})$, $(\tfrac{5\pi}{6}, \tfrac{5\pi}{6})$, 
$(\tfrac{\pi}{5}, \tfrac{2\pi}{5})$, $(\tfrac{2\pi}{5}, \tfrac{\pi}{5})$, $(\tfrac{3\pi}{5}, \tfrac{4\pi}{5})$, $(\tfrac{4\pi}{5}, \tfrac{3\pi}{5})$.

The correspondence between holomorphic data and points $(\ga,\de)$ is bijective if we fix $N\ (>0)$. For convenience we shall normalize by taking $N=1$.  Thus, the holomorphic data consists of $k_0,\dots,k_n$ with $0\le k_i+1\le 1$ and $\sum_{i=0}^{n} (k_i+1)=1$.  The 33 points satisfy $0\le k+1,l+1\le 1$, but only the 19 points with $a+b\le\pi$ satisfy 
$0\le k_i+1\le 1$ for all $i$, namely

(a) 15 points with $a,b\in\{ 0, \tfrac\pi3, \tfrac\pi2, \tfrac{2\pi}3, \pi \}$ and  $a+b\le\pi$;

(b) 4 additional points $(a,b)=(\tfrac{\pi}{6}, \tfrac{\pi}{6})$, $(\tfrac{\pi}{4}, \tfrac{\pi}{4})$, 
$(\tfrac{\pi}{5}, \tfrac{2\pi}{5})$, $(\tfrac{2\pi}{5}, \tfrac{\pi}{5})$.

\no These are the required 19 points.
\end{proof}

The five blocks in Table \ref{integralsolutions} divide the points into the following types (with reference to Fig.\ \ref{theregion}, Fig.\ \ref{thepoints}): top edge,  left hand edge,  diagonal edge, interior points on the central line of symmetry, then the remaining 4 interior points.
\begin{table}[h]
\renewcommand{\arraystretch}{1.5}
\begin{tabular}{c|c|c|c}
Cases 4a,4b & Cases 5a,5b & Cases 5c,5d,5e & Cases 6a,6b,6c
 \\
\hline
\sss $(3,1)$ & \sss $(4,2)$ & \sss $(3,1)$ & \sss $(4,2)$
\\
[-6pt]
\sss $(\tfrac{5}{3},1)$ & \sss $(\tfrac{7}{3},2)$ & \sss $(\tfrac{4}{3},1)$ & \sss $(2,2)$
\\
[-6pt]
\sss $(1,1)$ &  \sss $(\tfrac{3}{2},2)$ & \sss $(\tfrac{1}{2},1)$ &  \sss $(1,2)$
\\
[-6pt]
\sss $(\tfrac{1}{3},1)$ &  \sss $(\tfrac{2}{3},2)$ &  \sss $(-\tfrac{1}{3},1)$ &  \sss $(0,2)$
\\
[-6pt]
\sss $(-1,1)$ &  \sss $(-1,2)$ & \sss $(-2,1)$ &  \sss $(-2,2)$
\\
\hline
\sss $(-1,-\tfrac13)$ & \sss $(-1,\tfrac13)$ &  \sss $(-2,-\tfrac23)$ & \sss  $(-2,0)$
\\
[-6pt]
\sss  $(-1,-1)$ &  \sss $(-1,-\tfrac12)$ &  \sss $(-2,-\tfrac32)$ &  \sss $(-2,-1)$
 \\
 [-6pt]
\sss $(-1,-\tfrac53)$ &  \sss $(-1,-\tfrac43)$ &  \sss $(-2,-\tfrac73)$ &  \sss $(-2,-2)$
 \\
 [-6pt]
 \sss $(-1,-3)$ &  \sss $(-1,-3)$ &  \sss $(-2,-4)$ &  \sss $(-2,-4)$
 \\
\hline
\sss $(\tfrac{1}{3},-\tfrac{5}{3})$ & \sss $(\tfrac{2}{3},-\tfrac{4}{3})$ & \sss $(-\tfrac{1}{3},-\tfrac{7}{3})$ &  \sss $(0,-2)$
 \\
 [-6pt]
\sss $(1,-1)$ &  \sss $(\tfrac32,-\tfrac12)$ &  \sss $(\tfrac{1}{2},-\tfrac{3}{2})$ &  \sss $(1,-1)$
 \\
 [-6pt]
\sss $(\tfrac{5}{3},-\tfrac{1}{3})$ &  \sss $(\tfrac{7}{3},\tfrac{1}{3})$ &  \sss $(\tfrac{4}{3},-\tfrac{2}{3})$ &  \sss $(2,0)$
\\
\hline
\sss $(\tfrac13,-\tfrac13)$  &  \sss $(\tfrac23,\tfrac13)$ &  \sss $(-\tfrac{1}{3},-\tfrac{2}{3})$ &  \sss $(0,0)$
\\
[-6pt]
\sss $(0,0)$ &  \sss $(\tfrac14,\tfrac34)$ &  \sss $(-\tfrac34,-\tfrac14)$ &  \sss $(-\tfrac12,\tfrac12)$
\\
[-6pt]
\sss $(-\tfrac13,\tfrac13)$ & \sss $(-\tfrac16,\tfrac76)$ &  \sss $(-\tfrac76,\tfrac16)$  & 
 \sss $(-1,1)$
\\
\hline
 \sss $(1,-\tfrac13)$ &  \sss $(\tfrac32,\tfrac13)$ &  \sss $(\tfrac12,-\tfrac23)$ & \sss $(1,0)$
\\
[-6pt]
\sss $(\tfrac35,\tfrac15)$ & \sss $(1,1)$ &  \sss $(0,0)$ & \sss $(\tfrac25,\tfrac45)$
 \\
 [-6pt]
\sss  $(-\tfrac15,-\tfrac35)$ & \sss $(0,0)$ & \sss $(-1,-1)$ & \sss $(-\tfrac45,-\tfrac25)$
 \\
 [-6pt]
\sss $(\tfrac13,-1)$ & \sss $(\tfrac23,-\tfrac12)$ &  \sss $(-\tfrac13,-\tfrac32)$ & \sss $(0,-1)$
\end{tabular}
\bigskip
\caption{$(\ga,\de)$ for the 19 solutions with integral Stokes data.}\label{integralsolutions}
\end{table}

\newpage
\section{Holomorphic data and quantum cohomology}\label{interp}

The genus zero $3$-point Gromov-Witten invariants of a manifold $M$ lead to the quantum cohomology algebra $QH^\ast M$, and also to the quantum D-module $\calM$.  The latter is isomorphic to a D-module of the form $D^\la/I$, where $D^\la$ is a certain ring of differential operators and $I$ is a left ideal which depends on $M$.  It is equivalent to the Dubrovin/Givental connection.  We refer to \cite{CoKa99} or \cite{Gu08} for a detailed explanation of these concepts.

In many examples there is a natural presentation for the ideal $I$.
This is the case for the \ll small\rr  (orbifold) quantum cohomology of the variety $M=\X^{v_0,\dots,v_p}_{d_1,\dots,d_m}$ which is the intersection of hypersurfaces of degrees $d_1,\dots,d_m$ in weighted projective space $\P^{v_0,\dots,v_p}$.  It is known\footnote{Some assumptions on the hypersurfaces are necessary here. We refer to 
\cite{Be95},\cite{Gi96},\cite{CoGoXX},\cite{Go07},\cite{CCLT09}
for details.
See also \cite{GuSaXX}.
}
that $I$ is generated by a single differential operator.  This operator is obtained by left-dividing the operator
\[
\la^{\sum_0^p v_i}
\prod_{i=0}^p
v_i^{v_i} \b (\b-\tfrac1{v_i})\cdots(\b-\tfrac{v_i-1}{v_i}) -
\la^{\sum_1^m d_j}
\prod_{j=1}^m
d_j^{d_j} \b (\b-\tfrac1{d_j})\cdots(\b-\tfrac{d_j-1}{d_j})\, z
\]
by the highest common factor of the two summands.  Here, $\b=z\tfrac{d}{dz}$.  
In the quantum cohomology literature it is usual to write $z=q$, $\la=\hbar$, but we shall use $z,\la$ for consistency with earlier notation.

For example (Example 4.4 of \cite{GuSaXX}), in the case of the weighted projective space $\P^{1,2,3}$ itself, the $3$-point Gromov-Witten invariants determine and are determined by the Dubrovin/Givental connection
\[
\nabla=d+\frac1\lambda\,\eta\, dz=
d+\frac1\lambda
\begin{pmatrix}
  &  &  &  &  & \!\frac13z^{\frac13}
 \\
1  &   & & & &
 \\
  & 1 &  & & &
 \\
  & & \frac16 z^{\frac13} &  & &
 \\
  & & & \frac13 z^{\frac16} & &
 \\
  & & & & \frac12 z^{\frac16} & 
\end{pmatrix}
\frac{dz}z.
\]
The corresponding quantum D-module is defined by declaring that $\b$ acts on cohomology-valued functions as $\b+\frac 1\lambda\,z\eta$.  
This D-module is naturally isomorphic to
\[
D^\la/
\left(
2^23^3\la^6\b^3(\b-\tfrac13)(\b-\tfrac12)(\b-\tfrac23)-z
\right)
\]
because the identity element of the cohomology ring
is a cyclic element of the D-module and it is annihilated by 
the (action of the) operator $2^23^3\la^6\b^3(\b-\tfrac13)(\b-\tfrac12)(\b-\tfrac23)-z$.
Similarly, for a degree $2$ hypersurface $\X^{1,2,3}_2$ in $\P^{1,2,3}$, the quantum differential operator is the result of left-dividing
$
2^23^3\la^6\b^3(\b-\tfrac13)(\b-\tfrac12)(\b-\tfrac23)- 
2^2\la^2\b(\b-\tfrac12)z
$
by the common factor $2^2\la^2\b(\b-\tfrac12)$.  This gives 
$3^3\la^4\b^2(\b-\tfrac13)(\b-\tfrac23)- z$, which is in fact the
quantum differential operator of $\P^{1,3}=\X^{1,2,3}_2$. 

With this in mind, we shall express the holomorphic data corresponding to $(\ga,\de)$ as a differential operator of the form $\la^{n+1}T_k-z$, then consider whether this is a quantum differential operator of the above type.  

The operator $T_k$ is defined as follows for cases 4a, 4b (and the other cases are analogous).  First, let us write the holomorphic data as
\[
\frac1\la \, \eta \, {dz} =
\frac1\la
\bp
 & & & \!\!\!\! \!  z^{k_0+1} \\
 \!\! z^{k_1+1} & & & \\
  &  \!\!\!  z^{k_2+1} & & \\
  & &  \!\!\!  z^{k_3+1} &
\ep
\frac{dz}z.
\]
To calculate a corresponding scalar operator we must choose a cyclic element
of the D-module.  For our purposes, it will suffice to do this in following way (see
sections 4.2 and 6.3 of \cite{Gu08} for the general principles).  Let us write the
equations for parallel sections of the (dual) flat connection $d-\frac1\la \, \eta^t \, {dz} $ as
$\la\b Y^t = z\eta^t Y^t$, where
$Y=(y_0, y_1, y_2, y_3)$.  
We obtain four scalar equations
\[
z^{-(k_i+1)} \la\b\, z^{-(k_{i-1}+1)} \la\b\, z^{-(k_{i-2}+1)} \la\b\, z^{-(k_{i-3}+1)} \la\b\,y_i=y_i
\]
for $y_i$,  $i\in\Z$ mod $4$, and any of these four scalar operators would be suitable as $T_k$.
As a definite choice, we shall use
\[
T_k=\b(\b-(k_j\!+\!1))(\b-(k_j\!+\!k_{j+1}\!+\!2))(\b-(k_j\!+\!k_{j+1}\!+\!k_{j+2}\!+\!3))
\]
where  $k_j,k_{j+1},k_{j+2},k_{j+3}$
is (lexicographically) the lowest of the four possibilities.
Thus, we represent the D-module corresponding to the 
holomorphic data $k=(k_0,k_1,k_2,k_3)$ as
$D^\la/\left(  \la^{n+1} T_k - z \right)$.  

The operators $T_k$ are listed in
Tables \ref{theoremCi}-\ref{theoremCiiii} of the appendix.
For example, the solution labelled $(a,b)=(\tfrac\pi2,\tfrac\pi3)$ in case 4a has 
$k+1=k_0+1=\tfrac12$, $l+1=k_2+1=\tfrac 13$.  Since $k_1=k_3$ here
and $\sum_{i=0}^3(k_i+1)=1$, we have $k_1+1=k_3+1=\tfrac1{12}$.  
Hence 
\[
k+1=(k_0+1,k_1+1,k_2+1,k_3+1)=(\tfrac1{2},\tfrac1{12},\tfrac1{3},\tfrac1{12}).
\]
We choose 
$\tfrac1{12}, \tfrac1{3}, \tfrac1{12}, \tfrac1{2}$ as the lowest representative. This gives
$T_k=\b(\b-\tfrac1{12})(\b-\tfrac5{12})(\b-\tfrac6{12})$, as indicated in Table \ref{theoremCi}.

We now observe that {\em the holomorphic data of each solution on the top edge or left hand edge of the region of Fig. \ref{thepoints} can be interpreted as a quantum D-module $\calM$ of the above type.} The spaces $M$ are shown in Table \ref{completeintersections}.  

\begin{table}[h]
\renewcommand{\arraystretch}{1.5}
\begin{tabular}{c|c|c|c}
Cases 4a,4b & Cases 5a,5b & Cases 5c,5d,5e & Cases 6a,6b,6c
 \\
\hline
$\P^3=\P^{1,1,1,1}$ &  $\P^4=\P^{1,1,1,1,1}$ &  $\P^{1,1,1,2}$ & $\P^{1,1,1,1,2}$
\\
[0pt]
$\X^{1,1,1,6}_{2,3}$ &  $\X^{1,1,1,1,6}_{2,3}$ &    $\X^{1,1,6}_{3}$ &  $\X^{1,1,1,6}_{3}$
\\
[0pt]
$\X^{1,1,4}_{2}$ &  $\X^{1,1,1,4}_{2}$ &  $\P^{1,4}$ &  $\P^{1,1,4}$
\\
[0pt]
$\P^{1,3}$ &  $\P^{1,1,3}$ &  $\P^{2,3}$ &  $\P^{1,2,3}$
\\
[0pt]
$\P^{2,2}$ &  $\P^{1,2,2}$ &  $\P^{1,2,2}$ &  $\P^{2,2,2}$
\\
\hline
$\P^{1,3}$ &  $\P^{2,3}$ &  $\P^{1,1,3}$ &  $\P^{1,2,3}$
\\
[0pt]
$\X^{1,1,4}_{2}$ &  $\P^{1,4}$ &  $\X^{1,1,1,4}_{2}$ &  $\P^{1,1,4}$
 \\
 [0pt]
$\X^{1,1,1,6}_{2,3}$ &  $\X^{1,1,6}_{3}$ &  $\X^{1,1,1,1,6}_{2,3}$ &  $\X^{1,1,1,6}_{3}$
 \\
 [0pt]
$\P^3=\P^{1,1,1,1}$ &  $\P^{1,1,1,2}$ &  $\P^4=\P^{1,1,1,1,1}$ &  $\P^{1,1,1,1,2}$
\end{tabular}
\bigskip
\caption{Quantum cohomology interpretation for solutions with integral Stokes data.}\label{completeintersections}
\end{table}

For example,  the quantum differential operator of
$\X^{1,1,1,6}_{2,3}$ is obtained by left-dividing
\[
\la^9 6^6 \b^4(\b-\tfrac16)\cdots(\b-\tfrac56)-
\la^5 2^2 3^3 \b(\b-\tfrac12)\b(\b-\tfrac13)(\b-\tfrac23)z
\]
by  
$\la^5 \b^3(\b-\tfrac13)(\b-\tfrac12)(\b-\tfrac23)$. 
This gives the holomorphic data 
$T_k=\b^2(\b-\tfrac16)(\b-\tfrac56)$ for
the second solution in Table \ref{theoremCi}.  (We are ignoring the coefficients $6^6$, $2^23^3$; this corresponds to the normalization $c_0=\cdots=c_n=1$ of
the holomorphic data.)

Conversely, it can be verified that {\em every quantum differential operator 
for $\P^{v_0,\dots,v_p}$ or $\X^{v_0,\dots,v_p}_{d_1,\dots,d_m}$
of the form $\la^{n+1}T_k-z$ with order $4$, $5$, or $6$ appears in our tables.}   Let us state this more formally, as it gives a purely analytic characterization of certain quantum D-modules.
First, we remark that $D^\la/\left(  \la^{n+1} T_k - z \right)$ has the properties of an \ll abstract (orbifold) quantum D-module\rr  when $k$ satisfies the conditions

\no (Q) $k_i+1=0$ for at least one $i$,

\no (G) if $x$ belongs to $\{ k_j + 1, k_j + k_{j+1} + 2, \dots, k_j + \cdots + k_{j+n-1} + n \}$ then so does $1-x$. 

\no Property (Q) is motivated by $H^2 M\ne 0$ and property (G) by the grading of the orbifold quantum cohomology.  This generalizes the concept of abstract quantum D-module (as in Chapter 6 of \cite{Gu08}) to the orbifold case. It represents the \ll expected\rr local properties of a quantum D-module near $z=0$.  The quantum D-modules of the spaces 
$\P^{v_0,\dots,v_p}$, $\X^{v_0,\dots,v_p}_{d_1,\dots,d_m}$
certainly satisfy these conditions, but the converse is false, e.g. it is easy to see that
the abstract quantum D-module 
$D^\la/
\left(
\la^4\b^2(\b-\tfrac1{10})(\b-\tfrac9{10})-z
\right)$
does not arise from any $\P^{v_0,\dots,v_p}$ or $\X^{v_0,\dots,v_p}_{d_1,\dots,d_m}$.  

Thus the difficult question of characterizing the genuine quantum D-modules arises.  
In our --- admittedly very restricted --- situation, there is a simple answer: 
it follows from our calculations (Tables 
\ref{theoremCi}-\ref{theoremCiiii}) that they are
characterized by the property of having integral Stokes data:

\begin{corollary}\label{char}  For $n\in\{3,4,5\}$, assume that $k=(k_0,\dots,k_n)$ satisfies conditions {\em (Q)} and 
{\em (G)}.  Then: $D^\la/\left(  \la^{n+1} T_k - z \right)$ is isomorphic to the 
quantum D-module of a space of the form $\P^{v_0,\dots,v_p}$ or $\X^{v_0,\dots,v_p}_{d_1,\dots,d_m}$ if and only if $s_1^\R,s_2^\R$ are integers.
\end{corollary}

Regarding other solutions, we note that the case 
$
T_k=\b(\b-\tfrac1{n+2})(\b-\tfrac2{n+2})\cdots(\b-\tfrac {n}{n+2})
$
is associated to an unfolding of a singularity of type $A_n$. This case was considered in detail by Cecotti and Vafa.  For $n=4$ and $n=5$ these appear in Table \ref{theoremCi} and  \ref{theoremCii}/\ref{theoremCiii} 
respectively;  the solutions are interior points of Fig. \ref{thepoints}.

The trivial solution $u=v=0$ occurs in all cases, and corresponds to  $(\ga,\de)=(0,0)$,
$(s_1^\R,s_2^\R)=(0,0)$.  
It is always an interior point of the region 
(but not always on the central line of symmetry). 
After a change of variable of the form $z\mapsto z^p$, 
the holomorphic data for the trivial solution can be written in the form
\[
\bp
 & & & 1\\
 1 & & & \\
  & \ddots & & \\
   & & 1 & 
\ep
dz
\]

\newpage
\section{Appendix: tables of asymptotic,  Stokes, and holomorphic data}\label{app}

In Tables \ref{theoremCi}-\ref{theoremCiiii}, we list the 19 solutions of section \ref{integral}, indexed by $(a,b)$, together with the asymptotic data $(\ga,\de)$, the integral Stokes data $(s_1^\R,s_2^\R)$, and the holomorphic data $T_k$.  

For even dimensional matrices, the symmetry $(a,b)\mapsto (b,a)$ 
transforms $(\ga,\de)$ to $(-\de,-\ga)$
and preserves   $T_k$, so in Tables \ref{theoremCi}, \ref{theoremCiiii} we just list the 12 solutions with $a\ge b$, i.e.\ $\ga+\de\ge 0$.    

As in Table \ref{integralsolutions}, the five blocks in the tables group the points in this order: top edge,  left hand edge (omitted for even dimensional matrices),  diagonal edge, interior points on the central line of symmetry, other interior points.

\begin{table}[h]
\renewcommand{\arraystretch}{1.5}
\begin{tabular}{c||c|c|l}
$(a,b)=\pi(k\!+\!1,l\!+\!1)$ & $(\ga,\de)$ & $(s_1^\R,s_2^\R)$ & $\ \ \ \ \ \ \ \ T_k$
 \\
\hline
 $(\pi,0)$ &  $(3,1)$ &  $(\pm 4,-6)$ & $\b^4$
\\
 $(\tfrac{2\pi}{3},0)$ &  $(\tfrac{5}{3},1)$ &  $(\pm 3,-4)$ & $\b^2(\b-\tfrac16)(\b-\tfrac56)$
\\
 $(\tfrac{\pi}{2},0)$ &  $(1,1)$ &  $(\pm 2,-2)$ & $\b^2(\b-\tfrac14)(\b-\tfrac34)$
\\
 $(\tfrac{\pi}{3},0)$ &  $(\tfrac{1}{3},1)$ &  $(\pm 1,0)$ & $\b^2(\b-\tfrac13)(\b-\tfrac23)$
\\
 $(0,0)$ &  $(-1,1)$ &  $(0,2)$ & $\b^2(\b-\tfrac12)^2$
\\
\hline\hline
 $(\tfrac{\pi}{2},\tfrac{\pi}{2})$ &  $(1,-1)$ &  $(0,-2)$ & $\b^2(\b-\tfrac12)^2$
\\
 $(\tfrac{2\pi}{3},\tfrac{\pi}{3})$ &  $(\tfrac{5}{3},-\tfrac{1}{3})$ &  $(\pm 2,-3)$ & $\b^2(\b-\tfrac13)^2$
\\
\hline
 $(\tfrac{\pi}{3},\tfrac{\pi}{3})$ &  $(\tfrac13,-\tfrac13)$ &  $(0,-1)$ & $\b(\b-\tfrac16)(\b-\tfrac36)(\b-\tfrac46)$
\\
 $(\tfrac{\pi}{4},\tfrac{\pi}{4})$ &  $(0,0)$ &  $(0,0)$ & $\b(\b-\tfrac14)(\b-\tfrac24)(\b-\tfrac34)$
\\
 $(\tfrac{\pi}{6},\tfrac{\pi}{6})$ &  $(-\tfrac13,\tfrac13)$ &  $(0,1)$ & $\b(\b-\tfrac16)(\b-\tfrac36)(\b-\tfrac46)$
\\
\hline
 $(\tfrac{\pi}{2},\tfrac{\pi}{3})$ &  $(1,-\tfrac13)$ &  $(\pm 1,-2)$ & $\b(\b-\tfrac1{12})(\b-\tfrac5{12})(\b-\tfrac6{12})$
\\
 $(\tfrac{2\pi}{5},\tfrac{\pi}{5})$ &  $(\tfrac35,\tfrac15)$ &  $(\pm 1,-1)$ & $\b(\b-\tfrac15)(\b-\tfrac25)(\b-\tfrac35)$
\end{tabular}
\bigskip
\caption{Asymptotic, monodromy, and holomorphic data for cases 4a, 4b ($\ga+\de\ge 0$)}\label{theoremCi}
\end{table}

\newpage
\begin{table}[h]
\renewcommand{\arraystretch}{1.5}
\begin{tabular}{c||c|c|l}
$(a,b)=\pi(k\!+\!1,l\!+\!1)$ & $(\ga,\de)$ & $(s_1^\R,s_2^\R)$ & $\ \ \ \ \ \ \ \ T_k$
 \\
\hline
 $(\pi,0)$ &  $(4,2)$ &  $(5,-10)$ & $\b^5$
\\
 $(\tfrac{2\pi}{3},0)$ &  $(\tfrac{7}{3},2)$ &  $(4,-7)$ & $\b^3(\b-\tfrac16)(\b-\tfrac56)$
\\
 $(\tfrac{\pi}{2},0)$ &  $(\tfrac{3}{2},2)$ &  $(3,-4)$ & $\b^3(\b-\tfrac14)(\b-\tfrac34)$
\\
 $(\tfrac{\pi}{3},0)$ &  $(\tfrac{2}{3},2)$ &  $(2,-1)$ & $\b^3(\b-\tfrac13)(\b-\tfrac23)$
\\
 $(0,0)$ &  $(-1,2)$ &  $(1,2)$ & $\b^3(\b-\tfrac12)^2$
\\
\hline
 $(0,\tfrac{\pi}{3})$ &  $(-1,\tfrac13)$ &  $(0,1)$ & $\b^2(\b-\tfrac26)(\b-\tfrac36)(\b-\tfrac46)$
 \\
  $(0,\tfrac{\pi}{2})$ &  $(-1,-\tfrac12)$ &  $(-1,0)$ & $\b^2(\b-\tfrac14)(\b-\tfrac24)(\b-\tfrac34)$
 \\
  $(0,\tfrac{2\pi}{3})$ &  $(-1,-\tfrac43)$ &  $(-2,-1)$ & $\b^2(\b-\tfrac16)(\b-\tfrac36)(\b-\tfrac56)$
 \\
  $(0,\pi)$ &  $(-1,-3)$ &  $(-3,-2)$ & $\b^4(\b-\tfrac12)$
 \\
\hline
 $(\tfrac{\pi}{3},\tfrac{2\pi}{3})$ &  $(\tfrac{2}{3},-\tfrac{4}{3})$ &  $(-1,-1)$ & $\b^2(\b-\tfrac13)(\b-\tfrac23)^2$
 \\
 $(\tfrac{\pi}{2},\tfrac{\pi}{2})$ &  $(\tfrac32,-\tfrac12)$ &  $(1,-2)$ & $\b^2(\b-\tfrac14)(\b-\tfrac24)^2$
 \\
 $(\tfrac{2\pi}{3},\tfrac{\pi}{3})$ &  $(\tfrac{7}{3},\tfrac{1}{3})$ &  $(3,-5)$ & $\b^2(\b-\tfrac16)(\b-\tfrac26)^2$
\\
\hline
 $(\tfrac{\pi}{3},\tfrac{\pi}{3})$ &  $(\tfrac23,\tfrac13)$ &  $(1,-1)$ & $\b(\b-\tfrac16)(\b-\tfrac26)(\b-\tfrac36)(\b-\tfrac46)$
\\
 $(\tfrac{\pi}{4},\tfrac{\pi}{4})$ &  $(\tfrac14,\tfrac34)$ &  $(1,0)$ & $\b(\b-\tfrac18)(\b-\tfrac28)(\b-\tfrac48)(\b-\tfrac68)$
\\
 $(\tfrac{\pi}{6},\tfrac{\pi}{6})$ &  $(-\tfrac16,\tfrac76)$ &  $(1,1)$ & $\b(\b-\tfrac1{12})(\b-\tfrac2{12})(\b-\tfrac6{12})(\b-\tfrac8{12})$
\\
\hline
 $(\tfrac{\pi}{2},\tfrac{\pi}{3})$ &  $(\tfrac32,\tfrac13)$ &  $(2,-3)$ & $\b(\b-\tfrac1{12})(\b-\tfrac2{12})(\b-\tfrac3{12})(\b-\tfrac5{12})$
\\
 $(\tfrac{2\pi}{5},\tfrac{\pi}{5})$ &  $(1,1)$ &  $(2,-2)$ & $\b(\b-\tfrac1{10})(\b-\tfrac2{10})(\b-\tfrac4{10})(\b-\tfrac8{10})$
 \\
 $(\tfrac{\pi}{5},\tfrac{2\pi}{5})$ &  $(0,0)$ &  $(0,0)$ & $\b(\b-\tfrac15)(\b-\tfrac25)(\b-\tfrac35)(\b-\tfrac45)$
 \\
  $(\tfrac{\pi}{3},\tfrac{\pi}{2})$ &  $(\tfrac23,-\tfrac12)$ &  $(0,-1)$ & $\b(\b-\tfrac1{12})(\b-\tfrac4{12})(\b-\tfrac7{12})(\b-\tfrac8{12})$

\end{tabular}
\bigskip
\caption{Asymptotic, monodromy, and holomorphic data for cases 5a, 5b}\label{theoremCii}
\end{table}

\newpage
\begin{table}[h]
\renewcommand{\arraystretch}{1.5}
\begin{tabular}{c||c|c|l}
$(a,b)=\pi(k\!+\!1,l\!+\!1)$ & $(\ga,\de)$ & $(s_1^\R,s_2^\R)$ & $\ \ \ \ \ \ \ \ T_k$
 \\
\hline
 $(\pi,0)$ &  $(3,1)$ &  $(-3,-2)$ & $\b^4(\b-\tfrac12)$
\\
 $(\tfrac{2\pi}{3},0)$ &  $(\tfrac{4}{3},1)$ &  $(-2,-1)$ & $\b^2(\b-\tfrac16)(\b-\tfrac36)(\b-\tfrac56)$
\\
 $(\tfrac{\pi}{2},0)$ &  $(\tfrac{1}{2},1)$ &  $(-1,0)$ & $\b^2(\b-\tfrac14)(\b-\tfrac24)(\b-\tfrac34)$
\\
 $(\tfrac{\pi}{3},0)$ &  $(-\tfrac{1}{3},1)$ &  $(0,1)$ & $\b^2(\b-\tfrac26)(\b-\tfrac36)(\b-\tfrac46)$
\\
 $(0,0)$ &  $(-2,1)$ &  $(1,2)$ & $\b^3(\b-\tfrac12)^2$
\\
\hline
 $(0,\tfrac{\pi}{3})$ &  $(-2,-\tfrac23)$ &  $(2,-1)$ & $\b^3(\b-\tfrac13)(\b-\tfrac23)$
 \\
  $(0,\tfrac{\pi}{2})$ &  $(-2,-\tfrac32)$ &  $(3,-4)$ & $\b^3(\b-\tfrac14)(\b-\tfrac34)$
 \\
  $(0,\tfrac{2\pi}{3})$ &  $(-2,-\tfrac73)$ &  $(4,-7)$ & $\b^3(\b-\tfrac16)(\b-\tfrac56)$
 \\
  $(0,\pi)$ &  $(-2,-4)$ &  $(5,-10)$ & $\b^5$
 \\
\hline
 $(\tfrac{\pi}{3},\tfrac{2\pi}{3})$ &  $(-\tfrac{1}{3},-\tfrac{7}{3})$ &  $(3,-5)$ & $\b^2(\b-\tfrac16)(\b-\tfrac26)^2$
 \\
 $(\tfrac{\pi}{2},\tfrac{\pi}{2})$ &  $(\tfrac12,-\tfrac32)$ &  $(1,-2)$ & $\b^2
 (\b-\tfrac14)(\b-\tfrac24)^2$
 \\
 $(\tfrac{2\pi}{3},\tfrac{\pi}{3})$ &  $(\tfrac{4}{3},-\tfrac{2}{3})$ &  $(-1,-1)$ & $\b^2(\b-\tfrac13)(\b-\tfrac23)^2$
\\
\hline
 $(\tfrac{\pi}{3},\tfrac{\pi}{3})$ &  $(-\tfrac13,-\tfrac23)$ &  $(1,-1)$ & $\b(\b-\tfrac16)(\b-\tfrac26)(\b-\tfrac36)(\b-\tfrac46)$
\\
 $(\tfrac{\pi}{4},\tfrac{\pi}{4})$ &  $(-\tfrac34,-\tfrac14)$ &  $(1,0)$ & $\b(\b-\tfrac18)(\b-\tfrac28)(\b-\tfrac48)(\b-\tfrac68)$
\\
 $(\tfrac{\pi}{6},\tfrac{\pi}{6})$ &  $(-\tfrac76,\tfrac16)$ &  $(1,1)$ & $\b(\b-\tfrac1{12})(\b-\tfrac2{12})(\b-\tfrac6{12})(\b-\tfrac8{12})$
\\
\hline
 $(\tfrac{\pi}{2},\tfrac{\pi}{3})$ &  $(\tfrac12,-\tfrac23)$ &  $(0,-1)$ & $\b(\b-\tfrac1{12})(\b-\tfrac4{12})(\b-\tfrac7{12})(\b-\tfrac8{12})$
\\
 $(\tfrac{2\pi}{5},\tfrac{\pi}{5})$ &  $(0,0)$ &  $(0,0)$ & $\b(\b-\tfrac1{5})(\b-\tfrac2{5})(\b-\tfrac3{5})(\b-\tfrac4{5})$
 \\
 $(\tfrac{\pi}{5},\tfrac{2\pi}{5})$ &  $(-1,-1)$ &  $(2,-2)$ & $\b(\b-\tfrac1{10})(\b-\tfrac2{10})(\b-\tfrac4{10})(\b-\tfrac8{10})$
 \\
  $(\tfrac{\pi}{3},\tfrac{\pi}{2})$ &  $(-\tfrac13,-\tfrac32)$ &  $(2,-3)$ & $\b(\b-\tfrac1{12})(\b-\tfrac2{12})(\b-\tfrac3{12})(\b-\tfrac5{12})$
\end{tabular}
\bigskip
\caption{Asymptotic, monodromy, and holomorphic data for cases 5c,5d,5e}\label{theoremCiii}
\end{table}

\newpage
\begin{table}[h]
\renewcommand{\arraystretch}{1.5}
\begin{tabular}{c||c|c|l}
$(a,b)=\pi(k\!+\!1,l\!+\!1)$ & $(\ga,\de)$ & $(s_1^\R,s_2^\R)$ & $\ \ \ \ \ \ \ \ T_k$
 \\
\hline
 $(\pi,0)$ &  $(4,2)$ &  $(\pm 4,-5)$ &  $\b^5(\b-\tfrac12)$
\\
 $(\tfrac{2\pi}{3},0)$ &  $(2,2)$ &  $(\pm 3,-3)$ & $\b^3(\b-\tfrac16)(\b-\tfrac36)(\b-\tfrac56)$
\\
 $(\tfrac{\pi}{2},0)$ &  $(1,2)$ &  $(\pm 2,-1)$ & $\b^3(\b-\tfrac14)(\b-\tfrac24)(\b-\tfrac34)$
\\
 $(\tfrac{\pi}{3},0)$ &  $(0,2)$ &  $(\pm 1,1)$ & $\b^3(\b-\tfrac26)(\b-\tfrac36)(\b-\tfrac46)$
\\
 $(0,0)$ &  $(-2,2)$ &  $(0,3)$ & $\b^3(\b-\tfrac12)^3$
\\
\hline\hline
 $(\tfrac{\pi}{2},\tfrac{\pi}{2})$ &  $(1,-1)$ &  $(0,-1)$ & $\b^2(\b-\tfrac14)(\b-\tfrac24)^2(\b-\tfrac34)$
\\
 $(\tfrac{2\pi}{3},\tfrac{\pi}{3})$ &  $(2,0)$ &  $(\pm 2,-2)$ & $\b^2(\b-\tfrac16)(\b-\tfrac26)^2(\b-\tfrac46)$
\\
\hline
 $(\tfrac{\pi}{3},\tfrac{\pi}{3})$ &  $(0,0)$ &  $(0,0)$ & $\b(\b-\tfrac16)(\b-\tfrac26)(\b-\tfrac36)(\b-\tfrac46)(\b-\tfrac56)$
\\
 $(\tfrac{\pi}{4},\tfrac{\pi}{4})$ &  $(-\tfrac12,\tfrac12)$ &  $(0,1)$ & $\b(\b-\tfrac18)(\b-\tfrac28)(\b-\tfrac48)(\b-\tfrac58)(\b-\tfrac68)$
\\
 $(\tfrac{\pi}{6},\tfrac{\pi}{6})$ &  $(-1,1)$ &  $(0,2)$ & $\b(\b-\tfrac1{12})(\b-\tfrac2{12})(\b-\tfrac6{12})(\b-\tfrac7{12})(\b-\tfrac8{12})$
\\
\hline
 $(\tfrac{\pi}{2},\tfrac{\pi}{3})$ &  $(1,0)$ &  $(\pm 1,-1)$ & $\b(\b-\tfrac1{12})(\b-\tfrac3{12})(\b-\tfrac5{12})(\b-\tfrac6{12})(\b-\tfrac9{12})$
\\
 $(\tfrac{2\pi}{5},\tfrac{\pi}{5})$ &  $(\tfrac25,\tfrac45)$ &  $(\pm 1,0)$ & $\b(\b-\tfrac1{10})(\b-\tfrac2{10})(\b-\tfrac4{10})(\b-\tfrac6{10})(\b-\tfrac8{10})$
\end{tabular}
\bigskip
\caption{Asymptotic, monodromy, and holomorphic data for cases 6a, 6b, 6c ($\ga+\de\ge 0$)}\label{theoremCiiii}
\end{table}

{\em

\noindent
Department of Mathematics\newline
Faculty of Science and Engineering\newline
Waseda University\newline
3-4-1 Okubo, Shinjuku, Tokyo 169-8555\newline
JAPAN
   
   \noindent
Taida Institute for Mathematical Sciences\newline
Center for Advanced Study in Theoretical Sciences  \newline
National Taiwan University \newline
Taipei 10617\newline
TAIWAN

}

\end{document}